\documentclass[12pt,a4paper]{article}
\usepackage{amsmath,amsthm,amsfonts,amssymb,alltt}
\usepackage{graphics}
\usepackage{graphicx}
\usepackage{epsfig}

\newtheorem{theorem}{Theorem}[section]
\newtheorem{lemma}[theorem]{Lemma}

\theoremstyle{definition}

\newtheorem{example}[theorem]{Example}
\theoremstyle{remark}
\newtheorem{remark}[theorem]{Remark}

\numberwithin{equation}{section}

\begin{document}

\author{Vladimir Mityushev \\
Pedagogical University, \\ ul. Podchorazych 2, Krakow 30-084, Poland 
}
\title{Optimal packing of spheres in $\mathbb R^d$ and extremal effective conductivity}
\date{}
\maketitle

\begin{abstract}
Optimal packing of spheres in $\mathbb R^d$ is studied by optimization of the energy $E$ (effective conductivity) of composites with ideally conducting spherical inclusions. It is demonstrated that the minimum of $E$ over locations of spheres is attained at the optimal packing. The energy is estimated in the framework of structural approximations. This method yields upper bounds and sometimes exact values for the maximal concentrations of spheres in $\mathbb R^d$. A constructive algorithm for the optimal locations of spheres associated to the classes of the Delaunay graphs is constructed.    
\end{abstract}

Keyword: Optimal packing of spheres; discrete energy; effective conductivity, Voronoi tessellation; Delaunay graph

MSC: 52C17, 05B40, 74Q15

\section{Introduction}
Packing problems refers to geometrical optimization problems. Various methods from  different topics of mathematics were applied to these problems. Extended reviews can be found in the books \cite{Bezdek} \cite{Bezdek2}, \cite{CS}, \cite{Hales}, \cite{Toth}. One of the most popular problem is the optimal packing of spheres in $\mathbb R^d$ \cite{CS}. Its complete solution for 2D is given in \cite{Toth} and for 3D in \cite{Hales}.  

It was noted in \cite{Mit2012} that solution to the physical problem of the optimal effective conductivity in 2D implies solution to the geometrical problem of the packing disks. The physical problem can be stated as follows. Given ideal conductors (having infinite conductivity coefficient) of the shape $D_i \subset \mathbb R^d$ ($i=1,2,\ldots$). To locate $D_i$ ($i=1,2,\ldots$) in the host medium of a finite conductivity in such a way that the homogenized medium is macroscopically isotropic and its effective conductivity attains the minimal value. Rigorous mathematical statements of the physical problem with fixed geometry can be found in \cite{BP} and in other works devoted to homogenization.   

Recent results in structural approximations \cite{LB}, \cite{Kolpakov} shown that densely packed composites can be investigated by the functional associated to the discrete energy. It has the following structure
\begin{equation}
\min_{t_1, t_2, \ldots}\;\mathop{{\sum}}_{k,j} g_{kj}^{(0)} |t_k-t_j|^2,
\label{eq:en1}
\end{equation}
where the minimum is taken over the values $t_j$ prescribed to the inclusion $D_j$. The value  $g_{kj}^{(0)}$ expresses the main term of the interparticle flux between the neighbor inclusions $D_j$ and $D_k$ when the distance $\delta_{jk}$ between them tends to zero. Usually, $g_{kj}^{(0)}$ has a simple form. In the case of linear conductivity, $g_{kj}^{(0)}$ not always has a singularity as $\delta_{jk} \to 0$ (see 3D examples in \cite{Kolpakov}). For instance, the flux between two spheres in $\mathbb R^d$ for $d>3$ is always regular for linear conductivity. Fortunately, $g_{kj}^{(0)}$ is always singular for non-linear conductivity governed by the $p$-Laplacian for $p>\frac{d+1}2$. This fact enables us to consider the minimum \eqref{eq:en1} not only in $t_j$ but also in the locations of the spherical inclusions $D_j$. 

This paper is devoted to study the minimum \eqref{eq:en1} and its application to the sphere packing problem. Sec.\ref{sec2} shortly presents the structural approximation theory following \cite{Keller}, \cite{LB}, \cite{Kolpakov} in $\mathbb R^d$ for $d=2,3$. Sec.\ref{sec3} is devoted to extension of the theory to the general space $\mathbb R^d$ and construction of the corresponding discrete energy. Estimations of energy are performed in Sec.\ref{sec4}. Sec.\ref{sec5} contains concluding remarks and discussion.     

\section{Structural approximation in $\mathbb R^d$}
\label{sec2}
Let $\boldsymbol{\nu}_j$ ($1,2, \ldots, d$) be the fundamental translation vectors in the space $\mathbb R^d$ ($d\geq 2$) which form a lattice $\mathcal Q=\{\sum_{j=1}^d m_j \boldsymbol{\nu}_j:m_j \in \mathbb Z \}$. The fundamental parallelotope $Q_{0}$ is defined by its $2^d$ vertices $\frac 12 \sum_{j=1}^d (\pm\boldsymbol{\nu}_j)$. 
Tow points $\mathbf a, \mathbf b \in \mathbb R^d$ are identified if their difference $\mathbf a- \mathbf b=\sum_{j=1}^d m_j \boldsymbol{\nu}_j$ belongs to the lattice $\mathcal Q$. Hence, such a topology is introduced on $Q_{0}$ that the opposite faces are glued. In the case $\mathbb R^2$, the fundamental parallelogram $Q_{0}$ can be considered as the classical flat torus. We will call the introduced topology on $Q_{0}$ in $\mathbb R^d$ also by toroidal. The distance $\|\mathbf a-\mathbf b\|$ between two points $\mathbf a,\mathbf b \in Q_0$ is introduced as
\begin{equation}\|\mathbf a-\mathbf b\|:=\min_{m_1,\ldots,m_d \in \mathbb Z} \left|\mathbf a-\mathbf b+\sum_{j=1}^d m_j \boldsymbol{\nu}_j\right|,
\label{eq:sam1}
\end{equation} 
where the modulus means the Euclidean distance between the points $\mathbf a$ and $\mathbf b$ in $\mathbb R^d$.
\begin{figure}[htp]
\centering
\includegraphics[clip, trim=0mm 0mm 0mm 0mm, width=.7\textwidth]{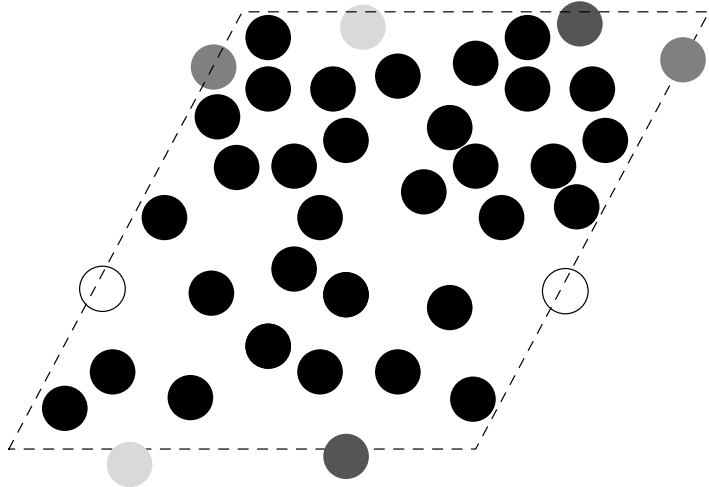} 
\caption{Periodicity cell $Q_{0}$. White and gray balls (disks) are identified in the torus topology.} 
\label{fig:Figure-1}
\end{figure}

Consider $n$ non-overlapping balls $D_k = \{\mathbf x \in \mathbb R^d: |\mathbf x - \mathbf a_k|<r \}$ of radius $r$ with the centers $\mathbf a_k$ in the cell $Q_{0}$ (see Fig.\ref{fig:Figure-1}). Let $D_0$ be the complement of all closure balls $D_k \cup \partial D_k$ to the domain $Q_{0}$. Following \cite{LB}, \cite{Kolpakov} and works cited therein we now shortly present the state of art in the non-linear homogenization of periodic composites governed by the  $p$-Laplace equation ($p\geq 2$)  
\begin{equation}
\nabla \cdot |\nabla u|^{p-2} \nabla u = 0, \quad \mathbf x \in D_{0}.  
\label{eq:p1}
\end{equation} 
The scalar function $-u(\mathbf x)$  is called the potential (sometimes the potential is taken as $u(\mathbf x)$), the vector function $\mathbf J=|\nabla u|^{p-2} \nabla u$ is called the flux. The $p$-Laplace equation \eqref{eq:p1} can be smoothly continued into the perforated domain $D_{0}+\sum_{j=1}^d \sum_{m_j\in \mathbb Z}  m_j \boldsymbol{\nu}_j$.  Introduce a unit vector $\boldsymbol{\xi} =(\xi_1,\xi_2,\ldots,\xi_d)\in \mathbb R^d$. The external potential is determined by the linear function $u_0(\mathbf x) =\boldsymbol{\xi} \cdot \mathbf x$ up to an arbitrary additive constant. The vector $\boldsymbol{\Omega}_{m_1,\ldots,m_d}=\sum_{j=1}^d m_j \xi_j \boldsymbol{\nu}_j$ with $m_j\in \mathbb Z$ corresponds to the potential jumps  with respect to the lattice $\mathcal Q$, i.e., the potential is quasi-periodic:
\begin{equation}
u\left(\mathbf x+\sum_{j=1}^d m_j \boldsymbol{\nu}_j\right)=u(\mathbf x)+\boldsymbol{\Omega}_{m_1,\ldots,m_d},\quad \mathbf x \in D_{0}\; (\forall m_j\in \mathbb Z).
\label{eq:p3}
\end{equation} 
The flux is periodic:
\begin{equation}
\mathbf J\left(\mathbf x+\sum_{j=1}^d m_j \boldsymbol{\nu}_j\right)=\mathbf J(\mathbf x), \quad \mathbf x \in D_{0}\; (\forall m_j\in \mathbb Z).
\label{eq:p3J}
\end{equation} 

The potential satisfies the boundary conditions
\begin{equation}
u(\mathbf x)=t_k, \quad |\mathbf x - \mathbf a_k|=r\; (k=1,2,\ldots,n),
\label{eq:p4}
\end{equation} 
where $t_k$ are undetermined constants. The total normal flux through each sphere vanishes:
\begin{equation}
\int_{\partial D_k}\mathbf J(\mathbf x)\cdot \mathbf n\; ds =0,  \quad  k=1,2,\ldots,n,
\label{eq:p5}
\end{equation}
where $\mathbf n$ denotes the outward unit normal vector to the sphere $\partial D_k$. 

The problem \eqref{eq:p3}-\eqref{eq:p5} describes the field in the periodic composite when the inclusions $D_k$ are occupied by perfect conductor and the host conductivity is governed by equation \eqref{eq:p1}. In electrostatics, $\mathbf E= \nabla u$ denotes the electric field, $\mathbf J= |\nabla u|^{p-2} \nabla u$ the electric current density. 
Energy passing through the cell per unit volume (effective conductivity) is calculated by formula (6.1.6) from \cite{LB}
\begin{equation}
\lambda=\frac 1{2|Q_0|} \int_{D_{0}} \mathbf J \cdot \mathbf E\; d\mathbf x  =\frac 1{2|Q_0|} \int_{D_{0}} |\nabla u|^{p} d\mathbf x,
\label{eq:p2}
\end{equation} 
where $|Q_0|$ stands for the volume of $Q_0$.
It is assumed that the centers $\mathbf a_k$ ($k = 1,2,\ldots,n$) are distributed in $Q_{0}$ in such a way that the corresponding composite is isotropic in macroscale, i.e., the effective conductivity of the composite is expressed by a scalar $\lambda$ and does not depend on the unit external flux $\boldsymbol{\xi}$. 

The energy $\lambda$ can be found as the minimum of the functional \cite{LB}: 
\begin{equation}
\lambda=\min_{v \in V}  \frac 1{2|Q_0|} \int_{D_{0}} |\nabla v|^{p} d\mathbf x,
\label{eq:p2e}
\end{equation} 
where the space $V$ consists of the quasi-periodic functions from the Sobolev space $W^{1,p}(Q_0)$:
\begin{equation}
V=\{v \in W^{1,p}(Q_0): v(\mathbf x)=t_k\; \mbox{on}\; D_k\; (k=1,2,\ldots,n)\}.
\label{eq:p2v}
\end{equation} 
Here, quasi-periodicity means that the conditions \eqref{eq:p3} and \eqref{eq:p3J} are fulfilled for $v$. Though the definition of the space $V$ depends on the external flux $\boldsymbol{\xi}$, the energy does not depend because of isotropy.    

The discrete network is a graph $\Gamma$ on the cell $Q_0$ with the vertices at $\mathbf a_k$ ($k=1,2,\ldots,n$) and the edges correspond to the necks  between neighbors. Neighbors are defined as  balls (disks) that share a common edge of the Voronoi tessellation of $Q_0$ with respect to their centers in toroidal  topology.  For each fixed $\mathbf a_k$, introduce the set $J_k$ of indexes for neighbor vertices and their total number $N_k=\# J_k$. 
The Voronoi tessellation in finite domains and the Delaunay graph $\Gamma$ are precisely described in \cite{LB}, \cite{Kolpakov}. The evident modification to periodic structures can be taken over. We use the term the Delaunay graph following \cite{LB}, \cite{Kolpakov} because it slightly differs from the Delaunay triangulation in degenerate cases. For example, consider a square and its four vertices. The traditional Delaunay triangulation has four sides of the square and one of the diagonals. In our approach, the Delaunay graph has only four sides (see Example \ref{ex4.3} in Sec.\ref{sec4}).

The discrete network model is based on the justification that the flux is concentrated in the necks between closely spaced inclusions (see \cite{LB} for non-linear equation \eqref{eq:p1} and \cite{LB}, \cite{Kolpakov} in the linear case $p=2$). First, we discuss the linear case when $p=2$. For two balls $D_k$ and $D_j$ the computation of the relative interparticle flux $g_{kj}$ (transport coefficient in terms of \cite{Kolpakov} and capacity in \cite{Keller}) relies on Keller's formulae \cite{Keller} 
\begin{equation}
g_{kj}=g_{kj}^{(0)} +O(\delta_{jk}^0), \quad \delta \to 0,
\label{eq:keller0}
\end{equation}
where in 3D
\begin{equation}
g_{kj}^{(0)}=-\pi r \ln \delta_{kj}
\label{eq:keller1}
\end{equation}
and in 2D
\begin{equation}
g_{kj}^{(0)}=\pi \sqrt{\frac{r}{\delta_{kj}}}.
\label{eq:keller2}
\end{equation}
Here, $\delta_{kj}$ denotes the gap between the balls $D_k$ and $D_j$:
\begin{equation}
\delta_{kj}=\|\mathbf a_k-\mathbf a_j\|-2r.
\label{eq:keller3}
\end{equation}
To each edge of the graph $\Gamma$ the flux \eqref{eq:keller1} or \eqref{eq:keller2} is assigned. 

In the non-linear case for $p>2$, we have the following formulae due to \cite{LB}
 in 3D
\begin{equation}
g_{kj}^{(0)}=\frac{\pi}{(p-2) r^{p-3}} \left(\frac{r}{\delta_{kj}} \right)^{p-2}
\label{eq:keller4}
\end{equation}
and in 2D
\begin{equation}
g_{kj}^{(0)}=\frac{(2p-5)!!}{(2p-4)!!} \;\frac{\pi^{\frac 32}}{r^{p-2}} \; \left(\frac{r}{\delta_{kj}} \right)^{p-\frac 32}.
\label{eq:keller5}
\end{equation}
It is worth noting that $g_{kj}^{(0)}$ depends only on geometry of the problem, i.e., on $r$ and on the gap $\delta_{kj}$ given by \eqref{eq:keller3}.  

\begin{remark}
The explicit coefficient $\frac{(2p-5)!!}{(2p-4)!!}$ is introduced in \eqref{eq:keller5} instead of the coefficient of the Taylor series for $(1-x)^{-\frac 12}$ in \cite{LB}. Formulae \eqref{eq:keller4} and \eqref{eq:keller5} are taken form \cite{LB} with slight corrections (see equations \eqref{eq:sa8a}, \eqref{eq:sa9a} and a remark below them).  
\end{remark}

\begin{figure}[htp]
\centering
\includegraphics[clip, trim=0mm 0mm 0mm 0mm, width=.7\textwidth]{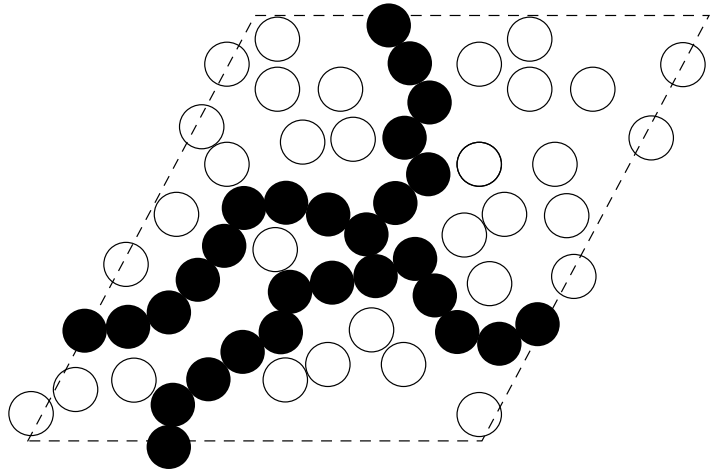} 
\caption{Percolation chains marked by black.} 
\label{fig:Figure-2}
\end{figure}

Let for shortness, $\mathbf t= (t_1,t_2, \ldots,t_n) \in \mathbb R^n$ and $\mathbf a = (\mathbf a_1,\mathbf a_2, \ldots, \mathbf a_n) \in \mathbb R^{n} \times \mathbb R^{d}$.
Introduce the double sum 
\begin{equation}
 \mathop{{\sum}'}_{k,j} = \sum_{k=1}^n \sum_{j\in J_k}.
\label{eq:sa2d}
\end{equation}
In order to formulate the main asymptotic result of \cite{LB} introduce the maximal length of edges $\delta = \max_k \max_{j \in J_k} \delta_{jk}$ of the graph $\Gamma$. Following \cite{LB}, \cite{Kolpakov} we consider the class $\mathcal D$ of macroscopically isotropic composites with densely packed inclusions. 
The term "densely packed inclusions" means that any such a location of balls has a percolation chain for $\delta=0$. More precisely, consider a set of non-overlapping balls in $Q_0$ with $\delta \geq 0$ endowed with the toroidal topology. Change of $\delta$ means that the centers of balls are fixed but their radii change and remain equal. Let for $\delta =0$ there exists a chain of touching balls connecting the opposite faces of the parallelotope $Q_0$. Such a chain is called a percolation chain\footnote{In the percolation theory, a percolation chain of balls is usually defined as a chain of overlapping open balls.}. Macroscopic isotropy implies that each pair of the opposite faces of $Q_0$ posses a percolation chain for $\delta =0$  (see Fig.\ref{fig:Figure-2}).
Any macroscopically isotropic location not belonging to the class $\mathcal D$ of densely packed inclusions can be replaced by an element of $\mathcal D$ having higher concentration. It is possible to do it by parallel translations of non-touching groups of balls to make them touched.

Introduce the main term of the discrete energy \cite{LB} 
\begin{equation}
\sigma =  \min_{\mathbf t}\;\frac 1{2|Q_0|} \mathop{{\sum}'}_{k,j} g_{kj}^{(0)} |t_k-t_j|^p,
\label{eq:sa1d}
\end{equation}
where $g_{kj}^{(0)}$ is given by \eqref{eq:keller4} or by \eqref{eq:keller5}. 

\begin{theorem}[\cite{LB}]
\label{thLB}
The continuous energy \eqref{eq:p2} and the discrete energy \eqref{eq:sa1d} for $p \geq 2$ in the spaces $\mathbb R^d$ ($d=2,3$)  tend to $+\infty$ as $\delta \to +0$. 
Moreover, $\lambda$ can be approximated by $\sigma$ for sufficiently small $\delta$: 
\begin{equation}
\lambda=\sigma +O(\delta^0), \;\mbox{as} \;  \delta \to 0.
\label{eq:sa2app}
\end{equation}
\end{theorem}

Th.\ref{thLB} justifies the theory of structural approximation \cite{LB}, \cite{Kolpakov} based on a mesoscopic discretization when edges and vertices of the graph $\Gamma$ correspond to inclusions. This theorem was applied in \cite{LB}, \cite{Kolpakov}, \cite{ryl2008} to systematic study of densely packed composites with given locations of inclusions.

\section{Discrete network in $\mathbb R^d$}
\label{sec3}
Though the theory of structural approximation \cite{LB}, \cite{Kolpakov}, \cite{ryl2008} was constructed and justified in $\mathbb R^2$ and in $\mathbb R^3$, its main results hold in the general space $\mathbb R^d$  under the condition that the  interparticle flux $g_{kj}$ tends to infinity as $\delta_{kj} \to 0$. It can be established by a formal repetition of the arguments from \cite{LB}, \cite{Kolpakov} in $\mathbb R^d$. In the present section, we partially fill this gap and calculate the main term $g_{kj}^{(0)}$ of $g_{kj}$ up to a constant order term $O(\delta_{jk}^0)$. 

Following Keller \cite{Keller} and explanations by Kolpakov \cite[p. 18-25]{Kolpakov} we proceed to estimate the interparticle flux $g_{jk}$ for two balls $D_j$ and $D_k$ with the centers located at the points $(0,0,\ldots,0,\pm \frac a2)$ of radius $r<\frac a2$ in $\mathbb R^d$ ($d =2,3, \ldots$). 
Let the potential $u(\mathbf x)$ is equal to constants $\pm C$ on the spheres $(x_d\mp \frac a2)^2+R^2=r^2$ where $R = \sqrt{x_1^2+x_2^2+\ldots +x_{d-1}^2}$. Near the gap the spheres can be approximated by the paraboloids \cite{LB}, \cite{Kolpakov}
\begin{equation}
x_n=\pm \left(\frac{\delta_{jk}}2+\frac{R^2}{2r}\right),
\label{eq:sph1a}
\end{equation}
where $\delta_{jk}=a-2r$.
The potential can be approximated up to $O(\delta_{jk}^0)$ by the function (see (6.4.12) from \cite{LB} for 2D and 3D cases)
\begin{equation}
u_0(\mathbf x) = \frac{x_d}{H(\mathbf x_{d-1})},
\label{eq:sa1a}
\end{equation}
where $H(\mathbf x_{d-1})=\delta +\frac{R^2}{r}$ is the distance between the paraboloids.
Following \cite{LB} we can approximate the local gradient by 
\begin{equation}
\nabla u_0(\mathbf x) =\left(0,\ldots,0, \frac{1}{\delta_{jk} +\frac{R^2}{r}}\right)+O(\delta_{jk}^0) 
\label{eq:saf1}
\end{equation}
and the local flux by
\begin{equation}
|\nabla u_0(\mathbf x)|^{p-2} \nabla u_0(\mathbf x) =\left(0,\ldots,0, \frac{1}{\left(\delta_{jk} +\frac{R^2}{r}\right)^{p-1}}\right)+O(\delta_{jk}^0) .
\label{eq:saf2}
\end{equation} 
Then, $g_{jk}$ is approximated up tp $O(\delta_{jk}^0)$ by the integral over the $(d-1)$-dimensional ball $B$ on the hyperplane $x_d=0$ defined by inequality $R<r$:
\begin{equation}
g_{jk}^{(0)} = \int_{B}\frac{d\mathbf x_{d-1}}{\left(\delta_{jk} +\frac{R^2}{r}\right)^{p-1}},
\label{eq:sa2a}
\end{equation}
where the differential $d\mathbf x_{d-1} = dx_1\;dx_2\; \cdots dx_{d-1}$. 

First, we consider partial cases of \eqref{eq:sa2a} discussed in the previous works \cite{Keller}, \cite{LB}, \cite{Kolpakov}. Let $p=d=2$. Then, \eqref{eq:sa2a} becomes
\begin{equation}
g_{jk}^{(0)} = \int_{-r}^r\frac{dR}{\delta +\frac{R^2}{r}} =\pi \sqrt{\frac{r}{\delta_{jk}}}. 
\label{eq:p2d2}
\end{equation}
Let $p=2$ and $d=3$. Then, \eqref{eq:sa2a} gives
\begin{equation}
g_{jk}^{(0)} = 2 \pi \int_0^r\frac{R\;dR}{\delta_{jk} +\frac{R^2}{r}} = \pi r \ln{\frac{r}{\delta_{jk}}}.
\label{eq:p2d3}
\end{equation}
Formulae \eqref{eq:p2d2}-\eqref{eq:p2d3} coincide with the corresponding formulae from \cite{Keller}, \cite{LB}, \cite{Kolpakov} (see also \eqref{eq:keller1}-\eqref{eq:keller2}).

We now proceed to investigate the general case. The integral \eqref{eq:sa2a} in the spherical coordinates $(R,\phi_1,\phi_2, \ldots, \phi_{d-2})$ becomes 
\begin{equation}
\begin{array}{lll}
g_{jk}^{(0)} = 
\\
\int_0^{2\pi}d\phi_{d-2} \int_0^{\pi}d\phi_{d-3} \cdots \int_0^{\pi}d\phi_{1} \int_0^{r} \frac{\sin^{d-3} \phi_1 \sin^{d-4} \phi_2 \cdots \sin \phi_{d-3}}{\left(\delta +\frac{R^2}{r}\right)^{p-1}}R^{d-2} dR.
\end{array}
\label{eq:sa3a}
\end{equation}
It can be calculated analogously to the volume of the $d$-dimensional ball. First, calculate
\begin{equation}
\int_0^{2\pi}d\phi_{d-2} \int_0^{\pi}d\phi_{d-3} \cdots \int_0^{\pi} \sin^{d-3} \phi_1 \sin^{d-4} \phi_2 \cdots \sin \phi_{d-3} \;d\phi_{1}=
\frac{2\pi^{\frac{d-1}2}}{\Gamma \left(\frac{d-1}2\right)},
\label{eq:si1}
\end{equation}
where the $\Gamma$-function is used. The integral in $R$ is calculated by formula 
\begin{equation}
\int_0^{r} \frac{R^{d-2} dR}{\left(\delta_{jk} +\frac{R^2}{r}\right)^{p-1}} = \frac{r^{d-1}\; 
{}_2 F_1\left(\frac{d-1}{2},p-1,\frac{d+1}{2},-\frac{r}{\delta_{jk}}\right)}{\delta_{jk} ^{ p-1}(d-1) },
\label{eq:si2}
\end{equation}
where the hypergeometric function ${}_2 F_1$ is used. Therefore,
\begin{equation}
g_{jk} = \frac{2\pi^{\frac{d-1}2} r^{d-1}}{\Gamma \left(\frac{d-1}2\right)} 
\frac{
{}_2 F_1\left(\frac{d-1}{2},p-1,\frac{d+1}{2},-\frac{r}{\delta_{jk}}\right)}
{\delta_{jk}^{ p-1}(d-1) } .
\label{eq:sa4a}
\end{equation}
Hereafter, we assume that $p$ is a natural number greater than 2 and 
\begin{equation}
p>\frac{d+1}2.
\label{eq:pd}
\end{equation} 
The following asymptotic formula takes place for $Z \to \infty$
\begin{equation}
{}_2F_1\left(\frac{d-1}{2},p-1,\frac{d+1}{2},-Z\right)=\frac{1}{Z^{\frac{d-1}{2}}}
\frac{\Gamma \left(\frac{d+1}2\right)\Gamma \left(p-\frac{d+1}2\right)}{\Gamma \left(p-1\right)}+O\left(\frac{1}{Z^{q}}\right),
\label{eq:sa5a}
\end{equation}
where $q>\frac{d-1}{2}$. Formula \eqref{eq:sa5a} is got by use of the package $Mathematica^{\circledR}$ Applying \eqref{eq:sa5a} to \eqref{eq:sa4a} we obtain the main asymptotic term of $g_{jk}$ as $\delta_{jk} \to 0$
\begin{equation}
g^{(0)}_{jk} =\frac{1}{\delta_{jk}^{p-\frac{d+1}{2}}} 
\frac{2 (\pi  r)^{\frac{d-1}{2}} 
\Gamma \left(\frac{d+1}{2}\right) 
\Gamma \left(p-\frac{d+1}{2}\right)}
{(d-1) \Gamma \left(\frac{d-1}{2}\right)
\Gamma (p-1)}.
\label{eq:sa6a}
\end{equation}
Let $d$ be an odd number. Then, \eqref{eq:sa6a} becomes
\begin{equation}
g^{(0)}_{jk} =\frac{1}{\delta_{jk}^{p-\frac{d+1}{2}}} 
\frac{(\pi  r)^{\frac{d-1}{2}}\left(p-\frac{d+3}{2}\right)!}
{(p-2)! }. 
\label{eq:sa8a}
\end{equation}
Let $d$ be an even number. Using formula
\begin{equation}
\Gamma \left(\frac{m}2\right)=\frac{(m-2)!! \sqrt{\pi}}{2^{\frac{m-1}2}}, \quad 
 m \in \mathbb N,
\label{eq:sa7a}
\end{equation} 
we rewrite \eqref{eq:sa6a} in the form\footnote{Our formulae \eqref{eq:sa6a}, \eqref{eq:sa8a} and \eqref{eq:sa9a} slightly diverge with the corresponding formulae from \cite{LB} (see for instance (6.2.13)) by multipliers. It can be related to the error power $p$ taken in \cite{LB} (e.g. (6.4.15) from \cite{LB} for 2D and 3D cases) instead of the correct $p-1$ from \eqref{eq:sa2a}.}

\begin{equation}
g^{(0)}_{jk} =\frac{1}{\delta_{jk}^{p-\frac{d+1}{2}}} 
\frac{\sqrt{\pi} (\pi  r)^{\frac{d-1}{2}}(2p-d-3)!!}
{2^{p-\frac{d+1}{2}}(p-2)}. 
\label{eq:sa9a}
\end{equation}

Introduce the functional associated with energy
\begin{equation}
E(\mathbf t, \mathbf a)= \frac 1{2|Q_0|} \mathop{{\sum}'}_{k,j} (t_k-t_j)^p
f(\|\mathbf a_k-\mathbf a_j\|)
\label{eq:sa1f}
\end{equation}
and the main term of the discrete energy
\begin{equation}
\sigma= \min_{\mathbf t}\; \frac 1{2|Q_0|} \mathop{{\sum}'}_{k,j} (t_k-t_j)^p
f(\|\mathbf a_k-\mathbf a_j\|).
\label{eq:sa1g}
\end{equation}
Here, the main term of the interparticle flux $g^{(0)}_{jk}$ is written through the function $f(x) = c (x-2r)^{-\left(p-\frac{d+1}{2}\right)}$ of one variable $x \geq 0$. The constant $c$ depends on $r$, $p$, $d$ and can be explicitly written by use of \eqref{eq:sa8a}, \eqref{eq:sa9a}.
Th.\ref{thLB} can be extended to the general space $\mathbb R^d$ as follows.

\begin{theorem}
\label{thLB2}
Consider the class $\mathcal D$ of densely packed balls in $\mathbb R^d$.
The continuous energy \eqref{eq:p2} (see also \eqref{eq:p2e}) and the discrete energy \eqref{eq:sa1g} for $p \in \mathbb N$ satisfying \eqref{eq:pd} tend to $+\infty$ as $\delta \to +0$. 
Moreover, $\lambda$ can be approximated by $\sigma$ for sufficiently small $\delta$: 
\begin{equation}
\lambda=\sigma +O(\delta^0), \;\mbox{as} \;  \delta \to 0.
\label{eq:sa1app}
\end{equation}
\end{theorem}
Proof of the theorem repeats the proof of Th.\ref{thLB} from \cite{LB} by its extension to the periodicity cell $Q_0$ in the space $\mathbb R^d$. 

\begin{theorem}
\label{th3}
Consider the class $\mathcal D$ of densely packed balls in $\mathbb R^d$.
Let $\sigma = \sigma(\mathbf a)$ attains the global minimum at a location $\mathbf a_*(\delta)$ for sufficiently small $\delta$. Then, the optimal packing is attained at $\mathbf a_*(0)$.
\end{theorem}

Proof. Let $\sigma(\mathbf a)$ attain the global minimum at $\mathbf a_*(\delta)$ for sufficiently small $\delta$ and $\phi_*$ denote the concentration of balls for the location $\mathbf a_*(0)$. The function $\sigma(\mathbf a)$ continuously depends on concentration \cite{LB}, \cite{Kolpakov}. Hence, the function $\sigma(\mathbf a_*(\delta))$ is continuous in $0<\delta<\delta_0$ for sufficiently small $\delta_0$ and $\sigma(\mathbf a_*(0))=+\infty$. 

Let the optimal packing be attained at another location $\mathbf a^*$ for which the concentration $\phi^*>\phi_*$. The location $\mathbf a^*$ contains a percolation chain. Take such a radius $r_0<r$ for which the concentration is reduced to $\phi_*$. Then, all the balls in this location $\mathbf a^*$ with the radius $r_0$ are separated from each other, hence the corresponding conductivity $\sigma(\mathbf a^*)$ is a finite number. But the minimal conductivity $\sigma(\mathbf a_*(\delta))$ tends to infinity as $\delta\to 0$. This yields a contradiction. 

The theorem is proved.

\section{Extremal energy}
\label{sec4}
Let all the periodic Delaunay graphs $\Gamma$ with $n$ vertices are divided onto the equivalence classes of isomorphic graphs. A class of graphs will be denoted by $\mathcal G$.  
We are looking for the global minimum of the functional \eqref{eq:sa1f} in $\mathbf t$ and $\mathbf a$ in the torus topology:
\begin{equation}
\min_{\mathbf t, \mathbf a} E(\mathbf t, \mathbf a)= \mathop{{\sum}'}_{k,j}
|t_k-t_j|^p f(\|\mathbf a_k-\mathbf a_j\|)
\label{eq:sa4}
\end{equation}  
in a fixed class $\mathcal G$.

It is evident that the minimum \eqref{eq:sa4} exists since the continuous function $f(x)$ decreases and $0 \leq f(\|\mathbf a_k-\mathbf a_j\|) \leq +\infty$  for all $\|\mathbf a_k-\mathbf a_j\| \geq 2r$.
The function $f(x)$ as a convex function for $2r \leq x <+\infty$ satisfies Jensen's inequality
\begin{equation}
\sum_{i=1}^M p_i f(x_i) \geq f \left(\sum_{i=1}^M p_i x_i \right),
\label{eq:sa5}
\end{equation} 
where the sum of positive numbers $p_i$ is equal to unity. Equality holds if and only if all $x_i$ are equal. Let the finite sum $\mathop{{\sum}'}_{k,j}$ in \eqref{eq:sa4}  is arranged in such a way that $x_i=\|a_k-a_j\|$ and $p_i=\frac 1T(t_k-t_j)^p$, where $T=\mathop{{\sum}'}_{k,j} (t_k-t_j)^p$. 

Application of \eqref{eq:sa5} yields
\begin{equation}
\mathop{{\sum}'}_{k,j} (t_k-t_j)^p f\left(\|\mathbf a_k-\mathbf a_j\|\right)
\geq T f \left(\frac 1T \mathop{{\sum}'}_{k,j} (t_k-t_j)^p \|\mathbf a_k-\mathbf a_j\| \right).
\label{eq:sa6}
\end{equation} 
H\"{o}lder's inequality states that for non-negative $a_i$ and $b_i$
\begin{equation}
\sum_{i=1}^M a_i b_i \leq  \left(\sum_{i=1}^M a_i^2 \right)^{\frac 12}  \left(\sum_{i=1}^M b_i^2 \right)^{\frac 12}.
\label{eq:sa7}
\end{equation} 
This implies that 
\begin{equation}
\mathop{{\sum}'}_{k,j} (t_k-t_j)^p \|\mathbf a_k-\mathbf a_j\| \leq \left[ \mathop{{\sum}'}_{k,j} (t_k-t_j)^{2p} \right]^{\frac 12} \left[ \mathop{{\sum}'}_{k,j} \|\mathbf a_k-\mathbf a_j\|^2 \right]^{\frac 12}.
\label{eq:sa8}
\end{equation} 
The function $f(x)$ decreases, hence \eqref{eq:sa6} and \eqref{eq:sa8} give
\begin{equation}
\mathop{{\sum}'}_{k,j} (t_k-t_j)^p f\left(\|a_k-a_j\|\right)
\geq T f \left(\frac 1T \left[ \mathop{{\sum}'}_{k,j} (t_k-t_j)^{2p} \right]^{\frac 12} \left[ \mathop{{\sum}'}_{k,j} \|a_k-a_j\|^2 \right]^{\frac 12} \right).
\label{eq:sa9}
\end{equation} 
The minimum of the right hand part of \eqref{eq:sa9} on $\mathbf a = (\mathbf a_1,\mathbf a_2, \ldots, \mathbf a_n)$ is achieved independently on $t_k$ for
\begin{equation}
\max_{\mathbf a} g(\mathbf a)=\mathop{{\sum}'}_{k,j} \|\mathbf a_k-\mathbf a_j\|^2.
\label{eq:sa10}
\end{equation}

\begin{lemma}
\label{lemm1}
For any fixed $n$, any local maximizer of $g(\mathbf a)$ is the global maximizer which fulfils the system of linear algebraic equations
\begin{equation}
\mathbf a_k=\frac 1{N_k} \sum_{j\in J_k} \mathbf a_j+\frac 1{N_k} \sum_{\ell=1}^d s_{j\ell}\boldsymbol{\nu}_{\ell},\quad k=1,2,\ldots,n,
\label{eq:sa11}
\end{equation}
where $s_{j\ell}$ can take the values $0,\pm 1$ in accordance with the fixed Voronoi tessellation (the fixed graph $\Gamma$). The system \eqref{eq:sa11} has always a unique solution $\mathbf a = (\mathbf a_1,\mathbf a_2, \ldots, \mathbf a_n)$ up to an arbitrary additive constant vector. 
\end{lemma}
 
Proof. It follows from the definition \eqref{eq:sam1} and the properties of the Voronoi tessellation that 
\begin{equation}
\|\mathbf a_k-\mathbf a_j\|=\left|\mathbf a_k-\mathbf a_j+\sum_{\ell=1}^d s_{j\ell}\boldsymbol{\nu}_{\ell}\right|
\label{eq:sam1a}
\end{equation} 
for some $s_{j\ell}$ which can take the values $0,\pm 1$.
The extremal points of \eqref{eq:sa10} can be found from the system of equations
\begin{equation}
\nabla_k g(\mathbf a)=0, 
\quad k=1,2,\ldots,n, 
\label{eq:sa12}
\end{equation}
where $\mathbf a_k= (x_1^{(k)}, x_2^{(k)}, \ldots,x_d^{(k)})$ and 
$$
\nabla_k = \left( \frac{\partial }{\partial x_1^{(k)}}, \frac{\partial }{\partial x_2^{(k)}}, \ldots, \frac{\partial }{\partial x_d^{(k)}}\right)
$$  
The parallelotope $Q_0$ is a compact manifold without boundary, hence all the extremal points of \eqref{eq:sa10} satisfy this system. Equations \eqref{eq:sa12} can be written in the equivalent form \eqref{eq:sa11}.     

One can see that the sum of all equations \eqref{eq:sa11} gives an identity, hence, they are linearly dependent. Moreover, if $\mathbf a = (\mathbf a_1,\mathbf a_2, \ldots, \mathbf a_n)$ is a solution of \eqref{eq:sa11}, then $(\mathbf a_1+\mathbf c,\mathbf a_2+\mathbf c, \ldots, \mathbf a_n+\mathbf c)$ is also a solution of \eqref{eq:sa11} for any $\mathbf c \in \mathbb R^d$. 
Consider the homogeneous system corresponding to \eqref{eq:sa11}
\begin{equation}
\mathbf a_k=\frac 1{N_k} \sum_{j\in J_k} \mathbf a_j, \quad
k=1,2,\ldots,n,
\label{eq:sa13}
\end{equation}
where $\mathbf a_k$ belong to the cell $Q_0$. The system \eqref{eq:sa13} can be decoupled by coordinates onto independent systems
\begin{equation}
x_p^{(k)}=\frac 1{N_k} \sum_{j\in J_k} x_p^{(j)},\; 
k=1,2,\ldots,n  \quad  (p=1,2,\ldots, d).  
\label{eq:sa14}
\end{equation}
Each $p$th system \eqref{eq:sa14} has only constant solutions
\begin{equation}
x_p^{(1)} = x_p^{(2)} = \cdots = x_p^{(n)}.
\label{eq:sa16}
\end{equation} 
This follows from the consideration of the quadratic form  
\begin{equation}
X_p=\frac 12 \mathop{{\sum}'}_{k,j} (x_p^{(k)}-x_p^{(j)})^2.
\label{eq:sa15}
\end{equation}
It is  symmetric and positive semi-definite. Therefore, the quadratic form \eqref{eq:sa14} has a global minimum attained at a linear set of $\mathbb R^n$. 
All the local minima of \eqref{eq:sa14} coincide with the global minimum. One can see that the quadratic form attains the global minimum at the set \eqref{eq:sa16}. Therefore, all solutions of \eqref{eq:sa14}, hence, of \eqref{eq:sa13} are only constants. Then, the system \eqref{eq:sa11} has only one condition of solvability for the right hand part which is fulfilled. Therefore, the system \eqref{eq:sa11} always has a unique solution up to an arbitrary additive constant.

The lemma is proved.

We now proceed to summarize the obtained results concerning the geometrical optimal packing problem of equal spheres stated  for periodic packing in $\mathbb R^d$.  We consider the periodicity cell $Q_0$ represented by a parallelotope with glued opposite faces which contains an arbitrary number $n$ of spheres. The centers of spheres $\mathbf a = \{\mathbf a_1,\mathbf a_2, \ldots, \mathbf a_n \}$ lie in $Q_0$, satisfy the non-ovelapping restriction $\|\mathbf a_k -\mathbf a_j\| \geq 2r$ and the corresponding Delaunay graph $\Gamma$ belong to a fixed class of graphs  $\mathcal G$. Lemma \ref{lemm1} gives a bound for the packing problem in the fixed class $\mathcal G$. This bound yields the exact global minimum when the inequality \eqref{eq:sa9} becomes  an equality. Let some of the differences $|t_k-t_j|$ vanish and the rest of the differences $|t_k-t_j|$ are equal. If for the equal differences $|t_k-t_j|$, the differences $\|\mathbf a_k -\mathbf a_j\|$ are also equal we have got this exact global minimum. Actually, such a situation frequently met, for instance for graphs corresponding to laminated lattices. 

\begin{example}
Consider the class of graphs $\mathcal A_2$ in $\mathbb R^2$ with $n=m^2$ ($m\in \mathbb N$) vertices containing the hexagonal lattice $A_2$. This class is determined by the kissing number $N_k=6$ for all the disks. The first step in the proof of optimality of the hexagonal lattice is the proof of equality $N_k=6$. Here, we refer to \cite{Toth}. Hence, if we restrict ourselves by the class $\mathcal A_2$, we do not go out a set of graphs containing the optimal packing. 

The second step consists in the direct check that the regular hexagonal lattice satisfies \eqref{eq:sa11}. After, a simple external flux has to be applied and the corresponding $|t_k-t_j|$ and $\|\mathbf a_k -\mathbf a_j\|$ have to be calculated. 
Consider the regular hexagonal lattice generated by the vectors 
\begin{equation} 
\boldsymbol{\nu}_1=\left(m,0 \right), \quad\boldsymbol{\nu}_2=\left(m\cos \frac{\pi}3, m\sin \frac{\pi}3 \right).
\label{eq:bas}
\end{equation} 
The radius of disks holds $r=\frac 12$. The periodicity cell $Q_0$ is determined by the vectors \eqref{eq:bas}. 
It is convenient to consider $A_2$ as a laminated lattice with layers perpendicular to the axis $x_2$. Introduce the number of the layer $q \in \mathbb Z$, where the $q$th layer consists of the disks whose centers have the $x_2$-coordinate equal to $\frac{\sqrt{3}}2 q$ (see Figure 3). 
\begin{figure}[htp]
\centering
\includegraphics[clip, trim=0mm 0mm 0mm 0mm, width=.7\textwidth]{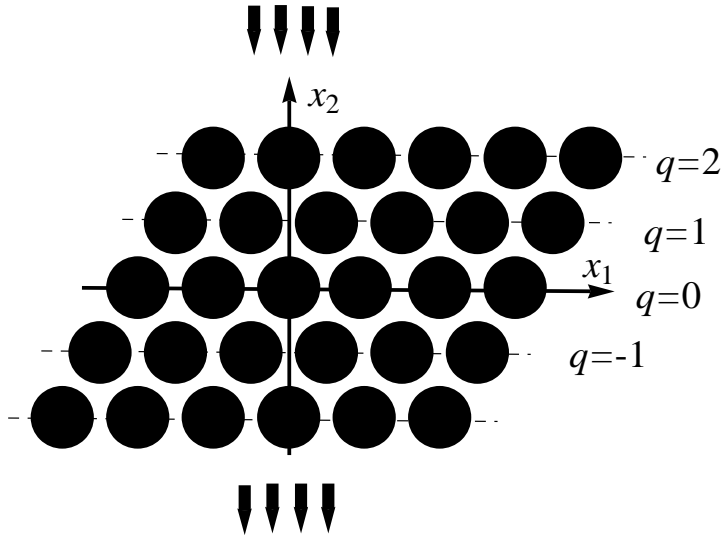} 
\caption{Hexagonal lattice as a laminated lattice with layers $q \in \mathbb Z$. The external flux is perpendicular to the axis $x_1$.} 
\label{fig:Figure-3}
\end{figure} 
Let the external flux be directed along the axis $x_2$ and it is determined by the external potential $u_0(\mathbf x)=\frac{\sqrt{3}}2 x_2$. Then the continuous and discrete potentials take the value $q$ on all the disks of the $q$th layer. Therefore, the difference $|t_k-t_j|$ is equal to zero if $\mathbf a_k$ and  $\mathbf a_j$ belong the same layer and it is equal to unity if $\mathbf a_k$ and  $\mathbf a_j$ belong neighbor layers. The difference $\|\mathbf a_k -\mathbf a_j\|$ takes the same value $2r=1$ for touching disks from neighbor layers. The inequality \eqref{eq:sa9} becomes an equality in this case. It follows from Lemma \ref{lemm1} that the global minimum \eqref{eq:sa4} on the class of graphs $\mathcal A_2$ is attains at the regular hexagonal graph. The optimal location $\mathbf a_*$ satisfies \eqref{eq:sa11}, hence depend on the basic vectors $\boldsymbol{\nu}_1$ and $\boldsymbol{\nu}_2$. If we take other basic vectors different from \eqref{eq:bas}, we arrive at the optimal location $\mathbf a_*$ which is not macroscopically isotropic. 

This example gives an alternative "physical" proof of the optimal 2D packing attained for the hexagonal array but the first step of the proof contains equation $N_k=6$ established by geometrical arguments. 
\end{example}

\begin{example}
\label{ex4.3}
Consider the class of graphs $\mathcal Z^2$ in $\mathbb R^2$ with $n=m^2$ ($m\in \mathbb N$) vertices containing the square lattice $\mathbb Z^2$. This class is determined by the contact number $N_k=4$ for all the disks. 

Consider the regular square lattice generated by the vectors $\boldsymbol{\nu}_1=\left(m,0 \right)$ and $\boldsymbol{\nu}_2=\left(0, m \right)$. The periodicity cell $Q_0$ is determined by these basic vectors. 
In accordance with the structural approximation theory the corresponding Delaunay graph\footnote{It does not formally produce a Delaunay triangulation in the commonly used sense.} consists of the edges parallel to the axes $x_1$ and $x_2$.   
The lattice $\mathbb R^2$ is laminated with layers perpendicular to the axis $x_1$. The $q$th layer consists of the disks whose centers have the $x_2$-coordinate equal to $q$. The external flux is determined by the external potential $-u_0(\mathbf x)= -x_2$. Then the continuous and discrete potentials take the value $q$ on all the disks of the $q$th layer. The further arguments repeat the previous example. Lemma \ref{lemm1} implies that the global minimum \eqref{eq:sa4} on the class of graphs $\mathcal Z^2$ is attains at the regular square graph.

It is worth noting that the regular square lattice $\mathbb Z^2$ gives the optimal packing in the class $\mathcal Z^2$. One can see that the class $\mathcal Z^2$ can be considered as unstable. It contains only one element $\mathbb Z^2$, since any perturbation of a vertex yields flipping and changes the structure of $\mathcal Z^2$. 
\end{example}

\begin{example}
Consider the class of graphs $\mathcal A_3$ in $\mathbb R^3$ with $n=m^3$ ($m\in \mathbb N$) vertices containing the regular face-centered cubic lattice $A_3$. Since this lattice has the laminated structure, we can consider the external flux directed perpendicular to the layers. Then, the potential is constant in each layer and the differences $\|\mathbf a_k -\mathbf a_j\| =2r$ are constant where the neighbor points $\mathbf a_k$, $\mathbf a_j$ belong to the neighbor layers. This implies that the lattice $A_3$ reaches the optimal packing in the class $\mathcal A_3$. It is worth noting that the optimal location $\mathbf a_*$ in the class $\mathcal A_3$ depends on the basic vectors prescribed to the fcc lattice $A_3$ since $\mathbf a_*$ satisfies the system \eqref{eq:sa11} whose the right hand side depends on the basic vectors. The class of graphs containing the lattice corresponding to the hexagonal close-packing determines another optimal location, of course the hcp. 

A class of graphs containing $A_3$ with arbitrary fixed basic vectors determines its optimal location $\mathbf a_*$ through the system \eqref{eq:sa11}. This location $\mathbf a_*$ depends on the basic vectors and must determine a macroscopically isotropic structure. The latter condition implies that the basic vectors have to correspond only to the fcc and hcp structures.  

Perhaps, the kissing number $N_k=12$ determines the class containing the optimal packing. This argument could give an alternative proof of the optimal packing for the regular fcc lattice \cite{Hales}. 
\end{example}

These examples can be extended to the general space. Consider the class of graphs $\mathcal A_d$ in $\mathbb R^d$ containing a regular laminated lattice $\Lambda_d$. Then, the regular lattice $\Lambda_d$ yields the optimal packing in the class $\mathcal A_d$. 

\section{Discussion} 
\label{sec5}
The most important fact used in this paper that the solution to the physical problem of the minimal energy (effective conductivity) implies solution to the geometrical problem of the optimal packing. Lemma \ref{lemm1} yields upper bounds for the classical optimal packing problem. For laminated structures these bounds are exact, hence solve the optimal packing problems within the considered classes.
Lemma \ref{lemm1} yields not only upper bounds but an effective algorithm for their computations, since the linear system \eqref{eq:sa11} is decomposed onto $d$ independent systems with $n$ equations.

Therefore, we solve the optimal packing problem in any fixed class of graphs $\mathcal G$. Let us discuss these classes of graphs. 
Though the number of spheres per cell $n$ is arbitrary and can tend to infinity, the periodicity cell is fixed at the beginning, hence the optimal location $\mathbf a_*$ depends on the basic vectors. It follows also from the observation that the right hand part of the system \eqref{eq:sa11} contains the basic vectors. For sufficiently large $n$, the majority of equations of \eqref{eq:sa11} are homogeneous and only "boundary" equations are inhomogeneous. Their number is of order $\sqrt{n}$. It is interesting to investigate the asymptotic dependence of $\mathbf a_*$ on the basic elements as $n \to \infty$. This will show, for instance, the effectiveness of the hexagonal type packing in the cubic cell in $\mathbb R^d$. 

We do not know a priory from the cell $Q_0$ and structure of $\mathcal G$ that the solution $\mathbf a_*$ of the system \eqref{eq:sa11}  will form a macroscopically isotropic structure. However, a simple necessary condition of isotropy can be checked: the location $\mathbf a_*$ must contain $d$ percolation chains connecting the opposite face of the parallelotope $Q_0$ (see Fig.\ref{fig:Figure-2}).

A straight forward verification which class $\mathcal G$ yields optimal packing can be proposed for a fixed $n$. First, for clarity consider 1D case. Then there exists exactly one class of graph $\mathcal Z$ including the regular 1D lattice $\mathbb Z$. A graph from $\mathcal Z$ consists of sequentially located edges along the axis. It follows from the system \eqref{eq:sa11} that each point $a_k$ lies in the middle of the segment $(a_{k-1}, a_{k+1})$ that immediately implies the points $a_k$ form the lattice $\mathbb Z$. The following extension of the above scheme to $\mathbb R^d$ can be proposed. Let $x_p^{(k)}$ denote the $p$th coordinate of the point $\mathbf a_k$ ($p=1,2,\ldots,d$). 
All the points $x_p^{(k)}$ ($k=1,2,\ldots,n$) are located on the real axis $\mathbb R$ (coincidence is permitted) and satisfy the system    
\begin{equation}
x_p^{(k)}=\frac 1{N_k} \sum_{j\in J_k} x_p^{(j)}+\frac 1{N_k} \sum_{\ell=1}^d s_{j\ell}\nu_p^{(\ell)},\quad 
k=1,2,\ldots,n,  
\label{eq:sa14b}
\end{equation}
where $\nu_p^{(\ell)}$ denotes the $p$th coordinate of the basic vector $\boldsymbol{\nu}_{\ell}$. The class of graphs $\mathcal G$ determines the sets $J_k$, i.e., connections between the points $x_p^{(k)}$ and $s_{j\ell}$. In particular, $\mathcal G$ determines $N_k$ restricted by the kissing number.   
It is important to solve the system \eqref{eq:sa14b} keeping the parameters $\nu_p^{(\ell)}$ in symbolic form. It can be done by separate numerical solutions to the systems
\begin{equation}
x_p^{(k \ell)}=\frac 1{N_k} \sum_{j\in J_k} x_p^{(j \ell)}+ \frac{s_{j\ell}}{N_k} ,\quad 
k=1,2,\ldots,n,   
\label{eq:sa14c}
\end{equation}  
and by the linear combination
\begin{equation}
x_p^{(k)}=\sum_{\ell=1}^d x_p^{(k \ell)}\nu_p^{(\ell)},\quad 
k=1,2,\ldots,n.  
\label{eq:sa14d}
\end{equation}  
All the computations can be made separately for $p=1,2,\ldots,d$. As a result we obtain expressions for the optimal $\mathbf a_k$ in the class $\mathcal G$   
\begin{equation}
\mathbf a_k=\sum_{\ell=1}^d \mathbf a_k^{(\ell)}\boldsymbol{\nu}_{\ell},\quad 
k=1,2,\ldots,n,  
\label{eq:sa14e}
\end{equation}  
where the vectors $\mathbf a_k^{(\ell)}$ are given numerically. The next numerical problem consists in determinations of such bases $\{\boldsymbol{\nu}_{\ell}\}_{\ell=1}^d$ which possess $d$ percolation chains connecting the opposite faces of the parallelotope $Q_0$. Ultimately, one has to perform this procedure for all the classes of graphs and take the best packing.     
The above algorithm is rather an idea how to apply solution to the packing problem within classes $\mathcal G$ to the traditional sphere packing \cite{CS}, \cite{Hales}.

The physical problem is considered in the class of periodic structures. General non-periodic composites were discussed in literature (see a review in \cite{Kolpakov} and references therein). The main question consists in extension of the structural approximation theory \cite{LB}, \cite{Kolpakov} developed for densely packed composites to the non-periodic case.  

We do not discuss the case of different radii. We suppose that a similar work can be done, since simple formulae for $g^{(0)}_{jk}$ are known \cite{LB}, \cite{Kolpakov} and can be extended to $\mathbb R^d$ following Sec.\ref{sec3}.







\end{document}